\documentclass[12pt,a4paper]{article}
\usepackage{amsfonts}
\usepackage{amsmath,amssymb,amsthm}
\usepackage[OT4]{fontenc}

\newtheorem{thm}{Theorem}[section]
\newtheorem{lem}[thm]{Lemma}
\newtheorem{prop}[thm]{Proposition}

\newtheorem{exam}[thm]{Example}

\newtheorem{f}[thm]{Fact}
\newtheorem{fs}[thm]{Facts}
\newtheorem{rem}[thm]{Remark}

\newtheorem{cor}[thm]{Corollary}
\newcommand{\To}{\stackrel{\circ}{\mathcal{T}}}

\newcommand{\mb}{\mathbb}
\newcommand{\mc}{\mathcal}

\newcommand{\gts}{\mathbf{GTS}}
\newcommand{\lss}{\mathbf{LSS}}
\newcommand{\sms}{\mathbf{Small}}
\newcommand{\ads}{\mathbf{ADS}}
\newcommand{\ds}{\mathbf{DS}}
\newcommand{\lds}{\mathbf{LDS}}
\newcommand{\wds}{\mathbf{WDS}}
\newcommand{\wss}{\mathbf{WSS}}
\newcommand{\spa}{\mathbf{Space}}
\newcommand{\vn}{\varnothing}

\begin{document}
\begin{center}

\Large
\textbf{On generalized topological spaces}

\vspace{3mm}
\large \textit{} \textsc{Artur Pi\k{e}kosz}
\end{center}

\begin{abstract}
In this paper a systematic study of the category $\gts$ of generalized topological spaces (in the sense of H. Delfs and M. Knebusch) and their   strictly continuous mappings begins. Some completeness and cocompleteness results are achieved. Generalized topological spaces help to reconstruct the important elements of the theory of locally definable and weakly definable spaces in the wide context of weakly topological structures. 
\end{abstract}
{\small \textit{2000 MS Classification:} 54A05, 18F10,  03C65.\\
\textit{Key words and frases:} generalized topological space, Grothendieck topology, topological structure.}

\section{Introduction}
Generalized topology in the sense od H. Delfs and M. Knebusch is an unknown chapter of general topology. In fact, it is a generalization of the classical concept of topology. The aim of this paper is to start the systematic study of generalized topology in this sense. (Do not mix with other meanings of "generalized topology" appearing in the literature.)

This concept was defined in \cite{LSS} and helped to develop a semialgebraic version of homotopy theory.
In 1991, in his  paper \cite{K91}, M. Knebusch suggested that his theory of locally semialgebraic spaces (developed in \cite{LSS} together with H. Delfs)   and weakly semialgebraic  spaces (developed in \cite{WSS}) 
could be generalized to the o-minimal context. (Here locally definable and weakly definable spaces allow to make constructions 
analogical to that known from the traditional homotopy theory.)

Successful generalization of the theory of Delfs and Knebusch 
to the case of o-minimal expansions of fields in \cite{ap}
opens a question if similar homotopy theories can be developed by the use of generalized topology and other ideas from \cite{LSS,WSS}.
Moreover, even if in many cases a full-fledged homotopy theory will not be available, the use of locally definable and weakly definable spaces 
over various structures will be important. And
even on the pure topological level, generalized topology is an interesting notion, worth to be studied.
Notice that the concept of generalized topology seems not to be  understood,
since the paper \cite{BO} suggests its non-existence (see \cite{ap} for a discussion) despite a heavy use of the theory of locally semialgebraic spaces from \cite{LSS}. The $\mc{T}$-spaces of \cite{EP} have natural generalized topologies, and it is convienient to see them as generalized topological spaces.

The category $\gts$ of generalized topological spaces and their strictly continuous mappings may be seen as an alternative to the usually used category \textbf{Top} of topological spaces and continuous mappings.
It originates from the categorical concept of Grothendieck topology, and contains \textbf{Top} as a full subcategory. In general, the third order concept od a generalized topological space is much more difficult to study than the second order concept of a topological space.
Only the use of $\gts$ allowed to glue infinitely many definable sets (considered with their natural o-minimal topologies) to produce locally definable spaces in \cite{LSS} and \cite{ap} and weakly definable spaces in \cite{WSS} and \cite{ap}. 
Many of the proofs of \cite{LSS} and \cite{WSS} are purely topological (naturally in the sense of the generalized topology), and it is important to extract the information on particular levels of structure.
The theory of infinite gluings of definable sets can be reconstructed to a large extent in the setting  of weakly topological structures. 

This paper is a continuation of the paper \cite{ap}, which was devoted to extending homotopy theory to the case of spaces over o-minimal expansions of fields. The present paper is much more general and gives a nice axiomatization and a basic theory of the category $\gts$. The new
(relative to traditional topology)
 concept of admissibility  is explained, and the paper deals with main (generalized) topological concepts as: small sets, bases, connectedness, completeness, paracompactness, Lindel\"ofness, separation axioms, and various concept of discreteness. 
We prove completeness and cocompleteness of the full subcategory $\sms$
of small spaces, and finite completeness of the subcategories $\lss$, $\wss_1$. 
Then $\gts$ is used to build some natural categories of spaces over model-theoretic structures. A natural setting here for topological considerations is assuming a topology on the underlying set $M$ of a model-theoretic structure $\mc{M}$, and insisting on regarding the product topologies on the cartesian powers of $M$. Such structures are called in this paper weakly topological. They are more general than so called first order topological structures of \cite{Pi,Ma} and structures with a definable topology of \cite{Sch}.
We prove cocompleteness of the category $\spa(\mc{M})$ of spaces over a weakly topological $\mc{M}$, and finite completeness of its full subcategories $\ads(\mc{M})$, $\ds(\mc{M})$, $\lds(\mc{M})$, $\wds_1(\mc{M})$.

The author hopes that from now on generalized topology (hidden in the language of locally semialgebraic spaces of \cite{LSS}, and weakly semialgebraic spaces of \cite{WSS}) will be developed without constraints.
Some open questions are suggested to the reader.
 
\vspace{2mm}
\textit{Notation.}
For families of sets $\mc{U},\mc{V}$, we will use usual operations $\cup,\cap,\setminus$ on sets, and analogical operations on families of sets, for example
$$\mc{U}\underline{\cap}\mc{V}=\{ U\cap V: U\in \mc{U}, V\in \mc{V}\},\:\:
\mc{U}\underline{\times}\mc{V}=\{ U\times V: U\in \mc{U}, V\in \mc{V}\}.$$
In particular, we will denote $V\underline{\cap} \:\mc{U}=\{ V\} \underline{\cap} \:\mc{U}$ for a set $V$.

\section{Grothendieck topology}

	Here we   remind what a Grothendieck topology is. The reader may consult books like \cite{bw,ks,mm}.

Let  $\mc{C}$ be a small category.
Consider the category of presheaves of sets on $\mc{C}$, denoted by $Psh(\mc{C})$ or $\hat{\mc{C}}$, which is the category of contravariant functors from $\mc{C}$ to $\mathbf{Sets}$.
Let us remind the fundamental fact:
\begin{f}[\textbf{Yoneda lemma, weak version}]
 Functor $\mc{C}\ni C \mapsto Hom(-,C)\in \hat{\mc{C}}$ is full and faithful, so we can consider $\mc{C}\subseteq \hat{\mc{C}}$. 
\end{f}

To define a Grothendieck topology the following notion is usually used:
a \textbf{sieve} $S$ on an object $C$ is a subobject of $Hom(-,C)$ as an object of $\hat{\mc{C}}$.
Since a sieve on $C$ is a presheaf of sets of morphisms with common codomain $C$, 
this may be translated (with a small abuse of language) into:  
$S$ is a set of morphisms with codomain $C$ such that $f\in S$ implies $f\circ g\in S$, if only $f\circ g$ is meaningful.
The largest sieve on $C$ is $Hom(-,C)$, the smallest is~$\emptyset$.  

A \textbf{Grothendieck topology} $J$ on $\mc{C}$ is a function $C\mapsto J(C)$, with $J(C)$ a set of sieves on $C$ such that the following axioms hold:

\textbf{(identity/nonemptyness)} for each $C$, $Hom(-,C)\in J(C)$;

\textbf{(stability/base change)} if $S\in J(C)$ and $f:D\rightarrow C$ is a morphism of $\mc{C}$, then
$f^*S=\{ g\mid f\circ g\in S\} \in J(D)$;

\textbf{(transitivity/local character)} if $S\in J(C)$ and $R$ is a sieve on $C$ such that $f^*R\in J(D)$ for each $f\in S$, $f:D\rightarrow C$, then $R\in J(C)$.

Elements of $J(C)$ are called \textbf{covering sieves}. Pair $(\mc{C},J)$ is called a \textbf{Gro\-then\-dieck site}.
The above axioms imply the conditions (see section III.2 of \cite{mm}):

\textbf{(saturation)} if $S\in J(C)$ and $R$ is a sieve containing $S$, then $R\in J(C)$;

\textbf{(intersection)} if $R,S \in J(C)$, then $R \cap S\in J(C)$.\\
It follows that  each $J(C)$ is a filter (not necessary proper) on the lattice $Sub_{\hat{\mc{C}}}Hom(-,C)$ of sieves on $C$.

If the category $\mc{C}$ \textit{has pullbacks}, then instead of covering sieves, we can speak about \textbf{covering families} of morphisms (generating respective sieves), so the axioms may be reformulated as:

\textbf{(identity/isomorphism)} for each $C$, $\{id_C \}$ is a covering family (stated also as: for each isomorphism
$f:D\rightarrow C$, $\{ f\}$ is a covering family);

\textbf{(base change)} if $\{ f_i:U_i\rightarrow U\}_i$ is a covering family, and $g:W\rightarrow U$ any morphism, then $\{ \pi_{2i}: U_i \times_U W \rightarrow W\}_i$ is a covering family;

\textbf{(local character)} if $\{ f_i:U_i\rightarrow U\}_i$ is a covering family, and $\{  g_{ij}:V_{ij}\rightarrow U_i\}_j$ are covering families, then $\{ f_i \circ g_{ij}:V_{ij}\rightarrow U\}_{ij}$ is a covering family;

and usually the following is added (see Definition 16.1.2 of \cite{ks}):

\textbf{(saturation)} if $\{ f_i:U_i\rightarrow U\}_i$ is a covering family, and
each of the $f_i$'s factorizes through an element of $\{ g_j:V_j\rightarrow U\}_j$ , then
$\{ g_j:V_j\rightarrow U\}_j$ is a covering family.

Alternatively (see section 6.7 of \cite{bw}) authors consider saturated and non-saturated Grothendieck topologies, but for a Grothendieck site saturation is usually assumed.

Grothendieck topology allows to define sheaves (of sets). A \textbf{sheaf} is such a presheaf $F$ that for each covering family $\{U_i\rightarrow U\}_i$ in the respective diagram
$$F(U) \stackrel{e}{\longrightarrow}\prod_i F(U_i)^{\stackrel{p_1}{\longrightarrow}}_{\stackrel{\longrightarrow}{p_2}} \prod_{i,j} F(U_i\times_U U_j)$$
the induced morphism $e$ is the equalizer of the standardly considered pair of
morphisms $p_1,p_2$ (cf. section III.4 in \cite{mm}).

A Grothendieck topology is \textbf{subcanonical} if every representable presheaf is a sheaf,
so in  this case we may, by identifying the objects of $\mc{C}$ with their respective representable presheaves,  consider $\mc{C}\subseteq Sh(\mc{C})\subseteq Psh(\mc{C})=\hat{\mc{C}}$. 
Grothendieck topologies used in practice are usually subcanonical.

\section{Generalized topological spaces}
This section is devoted to introducing a nice axiomatization and basic properties of generalized topological spaces.

For any set $X$, we have a boolean algebra $\mc{P}(X)$ of subsets of $X$, so this algebra of sets  may be treated as a small category with inclusions as morphisms. In this category
 fibered products are  the same  as (binary) products and the same as  (binary) intersections, so $\mc{P}(X)$ has pullbacks.
We want to introduce a full subcategory $Op$ of $\mc{P}(X)$, consisting of  ``open subsets'' of $X$.
Then we want to introduce a subcanonical Grothendieck topology on this category.
Subcanonicality means in this setting that for each covering family of morphisms
(which may be identified with a family of subsets of a given set, since morphisms are inclusions) the object covered by the family is the supremum of this family in the smaller category $Op$.

This leads to the notion of a generalized topological space introduced by H. Delfs and M. Knebusch in \cite{LSS}:
a \textbf{generalized topological space} (\emph{gts}) is a set $X$ together with a 
family of subsets $\stackrel{\circ}{\mathcal{T}} (X)$ of $X$,
called \textbf{open sets}, and a family of open families $\mathrm{Cov}_X $,
 called \textbf{admissible families} or
 \textbf{admissible (open) coverings},
such that:
\begin{description}
\item[(A1)] $\emptyset , X\in \To (X)$ 
(the empty set and the whole space are are open),

\item[(A2)] if $U_1 ,U_2 \in \To (X)$ then $U_1 \cup U_2 ,U_1
 \cap U_2 \in \To (X)$
(finite unions and finite intersections of open sets are open),

\item[(A3)] if $\{ U_i \}_{i\in I} \subset \To (X)$ and $I$ is finite, then
$\{ U_i \}_{i\in I} \in \mathrm{Cov}_X $
(finite families of open sets are admissible),

{The above three axioms are a strengthening of the identity axiom. They also insure that the smaller category has pullbacks. These three axioms may be collectively called the \textbf{finiteness} axiom.}

\item[(A4)] if $\{ U_i \}_{i\in I} \in \mathrm{Cov}_X$ then $\bigcup_{i\in I}
U_i \in \To (X)$ 
(the union of an admissible family is open),

{This axiom may be called \textbf{co-subcanonicality}. Together with subcanonicality, it means that admissible families are coverings (in the traditional sense) of their unions, which is imposed by the notation of \cite{LSS}:

$\{ U_i \}_{i\in I} \in \mathrm{Cov}_X (U)$ iff $U$ is the union of $\{ U_i \}_{i\in I}$.}

Subcanonicality and co-subcanonicality may be collectively called the \textbf{naturality} axiom.

\item[(A5)] if $\{ U_i \}_{i\in I} \in \mathrm{Cov}_X$, $V\subset
 \bigcup_{i\in I} U_i $, and $V\in \To (X)$, then $\{ V\cap  U_i \}_{i\in I}
 \in
 \mathrm{Cov}_X$
 (the intersections of an admissible family with an open subset of the union
 of the family form an admissible family),

{This is the stability axiom.}

\item[(A6)] if $\{ U_i \}_{i\in I} \in \mathrm{Cov}_X $ and for each $i\in I$
 there is $\{ V_{ij} \}_{j\in J_i} \in \mathrm{Cov}_X$ such that
 $\bigcup_{j\in J_i}
V_{ij} = U_i $, then $\{ V_{ij} \}_{\stackrel{i\in I}{ j\in J_i}}
 \in \mathrm{Cov}_X$ 
(all members of admissible coverings of members of an admissible
 family form together an admissible family),

{This is the transitivity axiom.}

\item[(A7)] if $\{ U_i \}_{i\in I} \subset \To (X)$, $\{ V_j \}_{j\in J} \in
 \mathrm{Cov}_X$, $\bigcup_{j\in J} V_j =\bigcup_{i\in I} U_i$,
 and  $\forall j\in J \: \exists i\in I : V_j
 \subset U_i$, then $\{ U_i \}_{i\in I} \in \mathrm{Cov}_X$
(a coarsening, with the same union, of an admissible family is admissible),

{This is the saturation axiom.}

\item[(A8)] if $\{ U_i \}_{i\in I} \in \mathrm{Cov}_X $, $V\subset
 \bigcup_{i\in I} U_i$ and $V\cap U_i \in \To (X)$ for each $i$, then
 $V\in \To (X) $
(if a subset of the union of an admissible family forms open
 intersections with the members of
 the family then the subset is open).
{This axiom may be called \textbf{regularity}.}

Both saturation and regularity have a smoothing character. Saturation may be achieved by adding coarsenings of admissible coverings for any \emph{gts}. Regularity may be achieved 
for the class of locally small spaces (see below for the definition) by adding ``locally open'' subsets and allowing locally essentially finite coverings (see section I.1 in \cite{LSS}). 
\end{description}

The above axioms may be restated shortly in the following way.
A \textbf{generalized topological space} is a triple
 $(M,Op_M,Cov_M)$, where $M$ is any set, $Op_M\subseteq \mc{P}(M)$,
and $Cov_M\subseteq \mc{P}(Op_M)$,  such that 
the following axioms are satisfied:

(\textbf{finiteness})\quad if $\mc{U}\in Fin (Op_M)$, then
$\bigcup \mc{U},\bigcap \mc{U}\in Op_M, \mc{U}\in Cov_M$,

(\textbf{stability})\ \quad\  if $V\in Op_M, \mc{U}\in Cov_M$, then $V\underline{\cap}\mc{U}\in Cov_{M}$, 

(\textbf{transitivity})\ if $\Phi \in \mc{P}(Cov_M), \underline{\bigcup} \Phi \in Cov_M$, then $\bigcup \Phi \in Cov_M$,  

(\textbf{saturation})\quad if $\mc{U}\in Cov_M, \mc{V}\in \mc{P}(Op_M), \mc{U}\preceq \mc{V}$, then $\mc{V}\in Cov_M$,

(\textbf{regularity})\quad  if $W\in \mc{P}(M),\mc{U}\in Cov_M$, $W\underline{\cap}\mc{U}\in \mc{P}(Op_M)$, then $W\cap(\bigcup\mc{U})\in Op_M$.     

\vspace{3mm}
In the above $Fin(\cdot)$ is the family of finite subsets of a given set, and $\mc{U}\preceq \mc{V}$ means: these two families have the same union and  $\mc{U}$ is a refinement of $\mc{V}$ (so $\mc{V}$ is a coarsening of $\mc{U}$).
  Notice that $Op_M = \bigcup Cov_M =\underline{\bigcup} Cov_M$, hence one can define a generalized topological space as just a pair $(M,Cov_M)$.

The naturality axiom does not appear, since our interpretation of the families $Op_M$ and $Cov_M$ is intended to give a Grothendieck site.
(Each member of $Cov_M$ covers its union.) 

A \textbf{strictly continuous mapping} between \emph{gts}es
 is such a mapping that the preimage of an admissible family is an admissible family. This, in particular, means that the preimage of an open set should be open.
(Strictly continuous mappings may be viewed as morphisms of sites, compare section 6.7 of \cite{bw} and section 17.2 of \cite{ks}.)
The \emph{gts}es together with the strictly continuous mappings form a category called here \textbf{GTS}. Isomorphisms of $\gts$ will be called \textbf{strict homeomorphisms}.

We will say that a family $\mc{U}$ is \textbf{essentially finite (countable)} if some finite (countable) subfamily $\mc{U}_0\subseteq \mc{U}$ covers the union of $\mc{U}$ (i. e. $\bigcup \mc{U}_0 =\bigcup \mc{U}$). 

\begin{exam}\label{przyk} We have the following simple examples of
\emph{gts}es:
\begin{enumerate}
\item
The space $\mb{R}^n_{alg}$, where  the closed sets are the algebraic subsets of $\mb{R}^n$, and the admissible coverings are the essentially finite open families (this is $\mb{R}^n$ with the Zariski topology).
\item
The space $\mb{C}^n_{alg}$, where  the closed sets are the algebraic subsets of $\mb{C}^n$, and the admissible coverings are the essentially finite open families (this is $\mb{C}^n$ with the Zariski topology).

\item
The space $\mb{R}^n_{salg}$, where the open sets are the open semialgebraic subsets of 
$\mb{R}^n$, and the admissible coverings are the essentially finite open families.

\item
The space $\mb{R}^n_{san}$, where  the open sets are the open semianalytic subsets of $\mb{R}^n$, and the admissible coverings are the open families essentially finite on bounded sets of $\mb{R}^n$.

\item
The space $\mb{R}^n_{suban}$, where  the open sets are the open subanalytic subsets of 
$\mb{R}^n$, and the admissible coverings are the open families
 essentially finite on bounded sets of $\mb{R}^n$.

\item
The space $\mb{R}^n_{top}$, the usual topological space $\mb{R}^n$.

\item The space $\mb{R}^n_{ts}$, where the open sets are the sets open in the usual topology, and the admissible coverings are the essentially finite open families.
\item
For each topological space $(X,\tau)$, we can take $Op_X=\tau$, and as $Cov_X$ the essentially countable open families.

\end{enumerate}
The function $\sin:\mb{R}\to\mb{R}$ is an endomorphism of $\mb{R}_{top}$ and of $\mb{R}_{san}$, but not an endomorphism of $\mb{R}_{salg}$.
\end{exam}

Understanding of the new concept of admissibility is given in the following two propositions and the next remark. First of all, notice that, for each open family $\mc{U}$, saturation implies: $\mc{U}\in Cov_M$ iff $\mc{U}\cup\{\vn\}\in Cov_M$ iff $\mc{U}\setminus \{\vn\}\in Cov_M$.

\begin{prop}
If $\mc{U}$ and $\mc{V}$ are admissible, then:

a) $\mc{U}\cup \mc{V}$ is admissible,

b) $\mc{U}\underline{\cup}\mc{V}$ is  admissible,
 
c) $\mc{U}\underline{\cap}\mc{V}$ is admissible.\\
Moreover, if $\bigcup\mc{U}$, $\bigcup \mc{V}$ are open, 
 $(\bigcup \mc{U})\cap (\bigcup \mc{V})=\emptyset$ and $\mc{U}\cup \mc{V}$ is admissible, then 

d) $\mc{U}$ and $\mc{V}$ are  admissible.
\end{prop}
\begin{proof}
a) Let $\Phi =\{ \mc{U}, \mc{V} \}$. Then, by finiteness, one has $\underline{\bigcup} \Phi =\{
\bigcup \mc{U}, \bigcup \mc{V}\}\in Fin(Op_M)\subseteq Cov_M$.
By transitivity, we get $\bigcup \Phi = \mc{U}\cup \mc{V}\in Cov_M$.

b) Notice that $\mc{U}\cup \mc{V}\preceq \mc{U}\underline{\cup}\mc{V}$ and use saturation.

c) Consider $\Phi=\{ \mc{U}\underline{\cap}V \mid V\in \mc{V}\}$.
Since $\underline{\bigcup} \Phi =\{ (\bigcup \mc{U})\cap V \mid V\in \mc{V}\}= (\bigcup \mc{U})\cap \mc{V}\in Cov_M$, we have
$\bigcup \Phi = \mc{U} \underline{\cap} \mc{V}\in Cov_M$.

d) This follows from stability and saturation, since $(\bigcup \mc{U} )\underline{\cap}(\mc{U}\cup \mc{V})\preceq \mc{U}$,
 and similarly for $\mc{V}$.
\end{proof}

\begin{prop}[omitting of admissible unions]
Assume the open families $\mc{U},\mc{V}_j (j\in J)$ are admissible and
the family $\mc{U}\cup \bigcup\limits_{j}(\mc{V}_j\cup \{\bigcup \mc{V}_j \})$ is admissible. Then $\mc{U}\cup \bigcup\limits_{j}\mc{V}_j$
is admissible.
\end{prop}
\begin{proof}
Notice that $\mc{U}\cup \bigcup\limits_{j}(\mc{V}_j\cup \{\bigcup \mc{V}_j \})$ is a refinement of $\mc{U}\cup \bigcup\limits_{j}\{ \bigcup \mc{V}_j\}$, so the second family is admissible. 
By applying transitivity to the family of families $\{\{U\} :U\in \mc{U}\} \cup \{ \mc{V}_j :j\in J\}$, we get admissibility of the family $\mc{U}\cup \bigcup\limits_{j}\mc{V}_j$.
\end{proof}

\begin{rem}
Notice that a subfamily of an admisssible family may not be admissible, even if they have the same union. 
Similarly if $\mc{U}$ and $\mc{V}$ are admissible, then $\mc{U}\cap\mc{V}$ may not be admissible.
For example, consider $\mb{R}$ to be an affine semialgebraic space 
(see \cite{ap}).
Take $\mc{U}=\{ (\frac{1}{n},1-\frac{1}{n}):n\geq 3 \}$, $\mc{V}=\{ (0,1)\} $, $\mc{W}=\{ (0,2)\}$.
Then $\mc{U}\cup \mc{V}$, $\mc{U}\cup \mc{W}\in Cov_{\mb{R}}$, but 
$\mc{U}=(\mc{U}\cup \mc{V})  \cap (\mc{U}\cup \mc{W})\notin Cov_{\mb{R}}$.

On the other hand, an open superfamily (with the same union) of an admissible family  is always admissible by saturation.
\end{rem}


We call a subset $K$ of a \emph{gts} $X$ \textbf{small}
 if for each admissible covering $\mc{U}$ of any open $U$, the set $K\cap U$
is covered by finitely many members of $\mc{U}$. 
 (We say in this case that the covering $\mc{U}$ is essentially finite on $K$ or on $K\cap U$.) The class  of small spaces forms a full subcategory \textbf{Small} of \textbf{GTS}.

\begin{prop} 
1) A subset of a small set is small. \\
2) The image of a small set by a strictly continuous mapping is small.
\end{prop}
\begin{proof}
1) Take any admissible (open) covering $\mc{U}$ of $U$. Assume $L\subseteq K$ and $K$ is small.
Since $\mc{U}$  is essentially finite on $K\cap U$, it is also essentially finite on 
its subset  $L\cap U$.

2) Assume $f:X\rightarrow Y$ is strictly continuous, and $K\subseteq X$ is small.
Take an admissible open covering $\mc{V}$  of $V$ in $Y$. Then $f^{-1}(\mc{V})$ is essentially finite on $K$, and $ f(K)\underline{\cap} \mc{V}= f(K\underline{\cap} f^{-1}(\mc{V}))$ is essentially finite.
\end{proof}

\begin{prop}
In any \emph{gts} if an  open family is  essentially finite (on its union), then this family is admissible.
\end{prop}
\begin{proof}
First notice that if an open family $\mc{U}$ has the largest element $U$, then
$\{ U\}$ and $\mc{U}$ are refinements of each other. Thus, by saturation, $\mc{U}$ is admissible.

Now if an open family $\mc{U}$ has a finite subcover $\mc{U}_0$ (of its union), then for each $V\in \mc{U}_0$
the family $V\underline{\cap} \mc{U}$ has the largest element $V$, and is admissible. Since $\mc{U}_0$ is admissible, also $\mc{U}_0\underline{\cap} \mc{U}$ is admissible by transitivity, and $\mc{U}$ is admissible by saturation.
\end{proof}

\begin{cor}
An open family of a small space is admissible if and only if it is essentially finite.
\end{cor}

\begin{exam}\label{e1}
Each topological space may be considered as a \emph{gts}. Admissible coverings are understood as any  coverings in the traditional sense.
Each continuous function is here strictly continuous. Thus \emph{\textbf{Top}} is a full subcategory of \emph{\textbf{GTS}}.
\end{exam}

Let us denote the small subsets of $X$ by $Sm_X$. We always have $Fin(X)\subseteq Sm_X\subseteq \mc{P}(X)$.
In the case of $X$ small, we get $Sm_X=\mc{P}(X)$. 
Denote $SmOp_X=Sm_X \cap Op_X$.



\begin{prop}\label{toprealline}
The usual euclidean space  with the usual topology (denoted $\mb{R}^n_{top}$) is a \emph{gts} in which all small subsets are finite. 
In particular, the compact interval $[0,1]$ in $\mb{R}_{top}$ is not small.
\end{prop}
\begin{proof}
A small subset $K$ of the (topological) space $\mb{R}^n_{top}$ is quasi-compact
and all of its subsets are also quasi-compact. But $\mb{R}^n_{top}$ is Hausdorff. Since all subsets of $K$ are compact, $K$ is a discrete set, and finally a finite set.
\end{proof}

\begin{exam}[non-examples]\label{e2}
Given a topological space, we could consider only singletons of open sets as admissible coverings of these sets. This would give so-called \textbf{indiscrete Grothendieck topology}, but \underline{not} (in general) a \emph{gts}, because of \emph{(A3)}. On the other hand, if we considered all families of open subsets as coverings of an open set, the resulting \textbf{discrete Grothendieck topology} would not  (in general) be subcanonical, hence \underline{not}  a \emph{gts}.
\end{exam}

Let us remind two interesting examples of \emph{gts}es.

\begin{exam}[``the subanalytic site,'' Remark 23 in \cite{ap}]\label{e3}
If $M$ is a real analytic manifold, then we can set:

a) an \emph{open subset} of a \emph{gts} means an open subanalytic subset;

b) an \emph{admissible covering} of a \emph{gts} means a family that is essentially finite on compact subsets.
\end{exam}
\begin{exam}[another ``site,'' Remark 23 in \cite{ap}]\label{e4}
In the situation as in Example \ref{e3}, set:

a) an \emph{open subset} of a \emph{gts} means any subanalytic subset;

b) an \emph{admissible covering} of a \emph{gts} means a family that is essentially finite on compact subsets.
\end{exam}

Notice that compact sets are small in Examples \ref{e3} and \ref{e4}, but not necessarily small in Example \ref{e1} (because an open subset of a compact set is not usually compact).

For any \emph{gts} $X$ and a small or open $Y\subseteq X$ we can induce a \emph{gts} on $Y$
(then $Y$ will be called a \textbf{subspace} of $X$)
 in the following way:

a) the open subsets of the \emph{gts} $Y$ are exactly the traces of open subsets of the \emph{gts} $X$ on $Y$,

b) the admissible families of the \emph{gts} $Y$ are exactly the traces of admissible families of the \emph{gts} $X$ on $Y$. 

This works for open subsets, since openness and admissibility in an open subspace is equivalent to openness and admissibility in the whole space. It works for small subsets because
in a small space  ``admissible'' means exactly ``essentially finite''. In general, the problem with transitivity, saturation and regularity arises.   

A subset $Y$ of a \emph{gts} $X$ is \textbf{closed} if its complement $Y^c$ is open, and is \textbf{locally closed} if it is an intersection of a closed set and an open set.  
A subset $Y$ of a \emph{gts} $X$ is \textbf{constructible} if it is a Boolean combination of open sets.
Remind that each constructible set is a finite union of locally closed sets.

If $\{(Op_{\alpha},Cov_{\alpha})\}_{\alpha \in A}$ is a family of
generalized topologies on a set $X$, then their intesection $(\bigcap_{\alpha} Op_{\alpha}, \bigcap_{\alpha} Cov_{\alpha})$
is a generalized topology. Thus, for any family $\{ V_{i}\}_{i\in I}$ 
of subsets of a \emph{gts} $X$ and any family $\{ \mc{U}_j \}_{j\in J}$
of subfamilies of $\{ V_{i}\}_{i\in I}$ (this is not restrictive because we may enlagre this family if needed),  we can speak about the generalized topology \textbf{generated by} $\{ V_{i}\}_{i\in I}$ and $\{ \mc{U}_j \}_{j\in J}$.
A \textbf{basis} of the generalized topology is such a family of open sets that any open set is an admissible union of elements from the basis.
(The notion of the basis of a topology is a special case of the notion of the basis of the generalized topology.)
Notice that if $\mc{B}$ is a basis of a \emph{gts} $(X,Op,Cov)$,
then $(\mc{B},Cov^{\mc{B}})$ generates $(Op,Cov)$, where
$Cov^{\mc{B}}= Cov \cap \mc{P}(\mc{B})$.

A subset $O\subseteq X$ will be called \textbf{weakly open} if it is a union of open subsets of the \emph{gts} $X$, and \textbf{weakly closed} if its complement is weakly open. The weakly open subsets form a topology generated by the open subsets.

The following example shows that a continuous mapping (in the generated topology) that maps small sets onto small sets may not be strictly continuous.

\begin{exam}
Take the real line $\mb{R}$ and define the open sets as the finite unions of open intervals
with  endpoints being rational numbers or infinities. Define the admissible coverings as the essentially finite coverings.
We get a small \emph{gts}. The mapping $\mb{R}\ni x\mapsto rx\in \mb{R}$, for $r\notin \mb{Q}$, is continuous but not strictly continuous.
\end{exam}

A subset $Y\subseteq X$ will be called \textbf{dense} in $X$ if its weak closure (i.e. closure in the generated topology) is equal to $X$. A \emph{gts} $X$ is \textbf{separable} if there is a countable dense subset of $X$.

We can introduce the countability axioms:
a \emph{gts} $X$ satisfies the \textbf{first axiom of countability} if each point of $X$ has a countable basis of open neighborhoods (such a family may easily be made admissible), and satisfies the \textbf{second axiom of countability} if there is a countable basis of the generalized topology. 

\begin{exam}
The semialgebraic real line $\mb{R}_{salg}$ has a countable basis of the 
generated topology $\{(p,q):p<q,\: p,q\in \mb{Q}\}$ which is not a basis of the generalized topology.
\end{exam}
Of course, each \emph{gts} satisfying the second axiom of countability is separable.

The separation axioms in \textbf{GTS} have weak and strong versions.
A \emph{gts} $X$ will be called:
a) \textbf{weakly $T_1$} if for each $x\in X$ and $y\in X\setminus\{x\}$ there is an open set $U$ such that $x\in U$ and $y\notin U$; b) \textbf{strongly $T_1$} if each singleton is a closed subset of $X$; c) \textbf{weakly Hausdorff} if for each pair of different points $x,y \in X$ there are open disjoint $U,V$ such that $x\in U$ and $y\in V$; d) \textbf{strongly Hausdorff}
if it is weakly Hausdorff and strongly $T_1$; e) \textbf{weakly regular} if for each point
$x\in X$ and each subset $F$ of $X$  not containing $x$ and being either closed or a singleton
there are open disjoint sets $U,V$ such that $x\in U$ and $F\subseteq V$; f) \textbf{strongly regular} if it is strongly $T_1$ and weakly regular; g) \textbf{weakly normal} if 
for two disjoint sets $F,G$ being each either closed or a singleton there are open disjoint sets $U,V$ such that $F\subseteq U$ and $G\subseteq V$;
h) \textbf{strongly normal} if it is strongly $T_1$ and weakly normal. 

It is clear that a \emph{gts} is weakly $T_1$ (weakly Hausdorff) if and only if the generated topology of $X$ is $T_1$ (Hausdorff, respectively).
Each weakly regular \emph{gts} has a regular generated topology.

\begin{exam}[strange small space]\label{wd}
Let $X$ be an infinite \emph{gts} where open sets are exactly finite sets or the whole space, and admissible coverings are exactly essentially finite open families.
Then $X$ is not strongly $T_1$, but the generated topology is discrete.
\end{exam}

\begin{prop} \label{gtsprod}
The category $\sms$ is complete.
\end{prop}
\begin{proof}
The product of a set of small spaces $\{X_{\alpha}\}_{\alpha \in A}$ is the cartesian product of sets $X=\prod_{\alpha \in A} X_{\alpha}$ with the generalized topology described as follows: open sets are finite unions of full cylinders with bases being finite cartesian products of open sets and the admissible coverings are the essentially finite coverings.
Then obviously the projections $\pi_{\alpha}:X\rightarrow X_{\alpha}$ are strictly continuous, and for a system of strictly continuous mappings $f_{\alpha}:Y\rightarrow X_{\alpha}$, the induced mapping $(f_{\alpha})_{\alpha}: Y\rightarrow X$ is strictly continuous.

The existence of equalizers for pairs of parallel strictly continuous functions between small spaces is obvious, since each subset of a small space forms a subspace. 
\end{proof}


There exist two inclusion functors $i_s : \sms \to \gts$ and $i_t: 
\mathbf{Top} \to \gts$, which are full and faithful.  
It is easy to check that the functor $i_s$ admits a left adjoint $sm$, which will be called the \textbf{smallification} functor:
$$sm(X)=X_{sm}, \quad sm(f)=f,$$
where $Op_{X_{sm}}=Op_{X},\quad Cov_{X_{sm}}=EssFin(Op_{X}) $.
Here $EssFin(\cdot)$ denotes the family of essentially finite families.

Similarly, the functor $i_t$ admits a right adjoint $top$ , which will be called the \textbf{topologization} functor:
$$top(X)=X_{top}, \quad top(f)=f,$$
where $Op_{X_{top}}=\tau(Op_{X}), \quad Cov_{X_{top}}=\mc{P}(\tau(Op_{X}))$.
Here $\tau(\cdot)$ denotes the generated topology.  
Hence $sm$ preserves colimits and $top $ preserves limits.
 
\begin{exam}[topologization of the boolean algebra of definable sets]
Let $\mc{M}$ be a first order mathematical structure (i. e. a set $M$ with some distinguished relations, constants and functions). 
Then for each $A\subseteq M$ and $n\in \mb{N}$, the family $Def_n(\mc{M},A)$ of $A$-definable subsets of $M^n$ forms a boolean algebra. Consider $Op=Def_n(\mc{M},A)$, and
$Cov=EssFin(Op)$. Then the generated topology $\tau(Op)$ is
the family of $\bigvee$-definable over $A$ subsets of $M^n$, and the complements of members of $\tau(Op)$ (i. e. the weakly closed sets) are the type definable over $A$ subsets of $M^n$.
(In practice, in model theory often bounds are set on the cardinality of the family of open sets forming a $\bigvee$-definable set.)
\end{exam}

\begin{f}[I.2 (p. 11) in \cite{LSS}]\label{gtscolim}
The category $\gts$ is cocomplete.
\end{f}

\begin{cor}
$\sms$ is a cocomplete category.
\end{cor}


An \textbf{open mapping} is a strictly continuous mapping between \emph{gts}es
such that the image of any open subset of the domain is an open subset of the range. 
It is clear that the canonical projections from the product of small spaces to its factors are open mappings. Similarly, a \textbf{closed mapping} maps closed subsets onto closed subsets.
The canonical projections from the product of small spaces to its factors are closed mappings.
Each strict homeomorphism is both an open mapping and a closed mapping.
A \textbf{local strict homeomorphism} is a strictly continuous mapping $f:X\to Y$ for which there is an open admissible covering $X={\bigcup}_{\alpha \in A} X_{\alpha}$ such that
each $f|_{X_{\alpha}}$ is an open mapping and a strict homeomorphism onto its image.
 The next example shows that a local strict homeomorphism may not be an open mapping.

\begin{exam}
A locally semialgebraic covering $p:\mb{R}_{loc}\to S^1$ (in the sense of  Example 13 in \cite{ap}) is a local strict homeomorphism but it is neither an open mapping nor a closed mapping. After passing to the strong topologies, we get the mapping $p_{top}:(\mb{R}_{loc})_{top}\to (S^1)_{top}$, which is open but not closed.
\end{exam}  


For any  family $\{X_{\alpha}\}_{\alpha \in A}$ of \emph{gts}es such that each intersection
$X_{\alpha}\cap X_{\beta}$ is an open subspace both in  $X_{\alpha}$ and in $X_{\beta}$, 
there is  a unique \emph{gts} $X$ having all $X_{\alpha}$s as open subspaces and the 
family   $\{X_{\alpha}\}_{\alpha \in A}$ admissible.
This is clear for finite families, and follows from the Fact \ref{gtscolim} for any families, since then $X$ is the colimit of the system of finite unions of spaces $X_{\alpha}$ partially ordered by inclusion.
In particular, if  an admissible covering $\mc{U}$ of $X$ is given, then the space $X$ is 
uniquely determined as the \textbf{admissible union} of the family $\mc{U}$ of open subspaces. In such a situation, we will write $X=\stackrel{a}{\bigcup} \mc{U}$. 
If the members of $\mc{U}=\{U_i\}_{i\in I}$ are pairwise disjoint, then the resulting space $X$ is caled the \textbf{direct (generalized topological) sum} (or coproduct) of the elements of $\mc{U}$, and will be denoted $X=\bigoplus_{i\in I} U_i$. Notice that then each $U_i$ is also closed. Moreover, each family of unions of $U_i$s is open (by regularity) and admissible (by saturation).

A \emph{gts} $X$ will be called \textbf{weakly discrete} if all its singlentons are open subsets,  \textbf{discrete} if all its subsets are open, and \textbf{topological discrete} if all families of subsets of $X$ are open and admissible. (For a topological discrete space $X$,
the formula $X=\bigoplus_{x\in X} \{x\}$ applies. All small topological discrete spaces are finite.) The space from Example \ref{wd} is weakly discrete but not discrete.

\begin{exam}
The semialgebraic space $\mb{R}$ has a small infinite weakly discrete subspace $\mb{N}$. The subspace $\mb{N}$ differs from the space from Example \ref{wd}. 
\end{exam}

\begin{exam}[infinite discrete small spaces]
On any infinite set $X$, there is still a generalized topology making $X$ a discrete small space. It is enough to set: open subset is any subset, admissible family is any essentially finite family.
\end{exam}

\begin{exam}
The identity mapping from the infinite discrete small space on some (infinite) set $X$ to the topological discrete space on this set is a closed and open  strictly continuous
bijection, but not a strict homeomorphism. 
\end{exam}

A \emph{gts} $X$ will be called \textbf{connected} if $X$ there is no  pair $U,V$ of open, disjoint, nonempty subsets of $X$ such that $U\cup V=X$.

\begin{prop}\label{cu}
If $X$ is a small space, $x\in X$, and $\{ C_{\alpha}\}_{\alpha\in A}$ is a family of connected subsets of $X$ each containing $x$, then $\bigcup_{\alpha \in A} C_{\alpha}$ is connected.
\end{prop}
\begin{proof}
Assume that the open subsets $U,V$ of $\bigcup_{\alpha \in A} C_{\alpha}$ are disjoint and nonempty, cover $\bigcup_{\alpha \in A} C_{\alpha}$, and $x\in U$. Let $y\in V$. Then
$y\in C_{\alpha_0}$ for some $\alpha_0 \in A$. The sets $U\cap C_{\alpha_0}$ and $V\cap C_{\alpha_0}$
are open, nonempty, and cover $C_{\alpha_0}$, thus $C_{\alpha_0}$ is not connected, contradiction. Hence $\bigcup_{\alpha \in A} C_{\alpha}$ is connected.
\end{proof}

The \textbf{connected component} of a point $x\in X$ of a small space $X$ is
the largest connected set $C_x$ containing $x$. Since the weak closure of $C_x$ is connected, 
$C_x$ is weakly closed. The \textbf{quasi-component} $\tilde{C}_x$ of $x$ is the intersection of clopen subsets of $X$ containing $x$. Each quasi-component is also weakly closed.
Each connected component is contained in a quasi-component.

We will say that a space $X$ \textbf{satisfies} $(\mathbf{ACC})$ if the family of connected components
of $X$ is open and admissible.
\begin{prop}
Each small space satisfying $(ACC)$ has a finite number of connected components.
\end{prop}
\begin{proof}
Choose one point in any connected component of the space. The resulting subspace is topological discrete and small, so it is finite.
\end{proof}



\section{Locally small spaces}
In this chapter we rebuild the theory of locally semialgebraic spaces from
\cite{LSS} on a pure topological level.

A \emph{gts} is \textbf{locally small} if there is an admissible covering of the whole space by  small open subsets. 
In other words: a \emph{gts} $X$ is locally small iff $X=\stackrel{a}{\bigcup} SmOp_X$.
Locally small \emph{gts}es form a full subcategory $\lss$ of $\gts$. 
Notice that in a locally small \emph{gts} each small set is contained in a small open set.
A family of subsets of a locally small space is called \textbf{locally finite} if each open small set meets only finitely many members of the family.

\begin{exam}
Each topological space $(X, \tau)$ can be considered as a locally small space relative to a chosen open covering $\mc{U}$ of the whole of $X$ in the following way: we declare members of $\mc{U}$ to be small spaces, and
the generalized topology on $X$ is then given by the formula $X=\stackrel{a}{\bigcup} \mc{U}$.
In particular, each topological discrete space is locally small (and it is not small, if infinite).
\end{exam}

Notice that $\mc{T}$-spaces of \cite{EP} are (practically)
locally small spaces in the above sense, which is the purpose of
introducing the family $\mc{T}_{loc}$ and the Grothendieck site $X_{\mc{T}_{loc}}$ in \cite{EP}.

For small spaces of type $\mb{R}^n_{sth}$, with the underlying set $\mb{R}^n$, we will consider their \textbf{localization} in the following sense:
assume that each of the open balls $B_n$ centered at the origin with radius $n\in\mb{N}\setminus \{0\}$ is open, and consider the localization $(\mb{R}^n_{sth})_{loc}$ to be the admissible union of this family of open small balls. 

\begin{exam}
The spaces $\mb{R}_{san}^n$, $\mb{R}_{suban}^n$ from Example \ref{przyk} as well as the localizations $(\mb{R}_{salg}^n)_{loc}$, $(\mb{R}^n_{ts})_{loc}$ are locally small, but not small.
\end{exam}

\begin{prop}
An open  family of a locally small space is  admissible if and only if it is ``\textbf{locally essentially finite}'',  which means:  the family is essentially finite on each small open subset.
\end{prop}
\begin{proof}
Each admissible family is locally essentially finite by the definition of a small subset.
If an open family $\mc{V}$  is locally essentially finite, then it is essentially finite on members of an admissible covering $\mc{U}$ of the space by open small subspaces. It means that for each member $U$ of   $\mc{U}$, the family $U\underline{\cap} \mc{V}$ is admissible. By the transitivity axiom, the family $\mc{U}\underline{\cap} \mc{V}$ is admissible. By the saturation axiom, the family $\mc{V}$ is admissible.   
\end{proof}

\begin{cor}
Each locally finite open family in a locally small space is admissible.
\end{cor}

\begin{prop}\label{closed}
Each locally essentially finite union of closed subsets of a locally small space is closed.
\end{prop}
\begin{proof}
Let $Z$ be a locally essentially finite union of closed sets. Take an admissible covering $\mc{U}$ of the space by small open subsets. For each element $U$ of $\mc{U}$, the set $Z\cap U$
is a finite union of relatively closed subsets of $U$, so it is relatively closed.
Each relative  complement $Z^c \cap U$ is open in $U$. By regularity, the set $Z^c$ is open. 
\end{proof}

\begin{lem}\label{lokmal}
The category $\lss$ has finite products.
\end{lem}
\begin{proof}
For a family $X_1,...,X_k$ of locally small spaces assume that their admissible coverigs
by small open subsets $\mc{U}_1,...,\mc{U}_k$, respectively,  are given.
Then the product space is given by the formula
$$X_1\times...\times X_k = \stackrel{a}{\bigcup} \mc{U}_1\underline{\times}...\underline{\times} \mc{U}_k.$$ 
The projections $\pi_i:X_1\times ...\times X_k \to X_i$ are obviously strictly continuous, and for given strictly continuous $f_i:Y\to X_i$, the mapping $(f_1,...,f_k):Y\to X_1\times ...\times X_k$ is also strictly continuous, since ``admissible'' means ``locally essentially finite''.
\end{proof}

The next example shows that the canonical projections along locally small spaces may not be open.
\begin{exam}
Consider the projection $\pi_2: \mb{R}_{loc}\times \mb{R}\rightarrow \mb{R}$ (as locally semialgebraic spaces, see \cite{LSS} or \cite{ap}).
The image of an open locally semialgebraic set $\bigcup_{n\in \mb{N}} (n,n+1)\times(n,n+1)$
under $\pi_2$ is not an open subspace of $\mb{R}$.
\end{exam} 
\begin{f}
If $\mc{U}=\{U_i\}_{i\in I}$ is a locally finite open family in some locally small space, and for each $i\in I$
some open $V_i\subseteq U_i$ is given, then $\{ V_i\}_{i\in I}$ is locally finite, thus admissible.
\end{f}

\begin{f}
If $X$ is a small space and $Y$ is a locally small space, then
the projection $X\times Y\rightarrow  Y$ is an open and closed mapping.
\end{f}

On locally small spaces we introduce a topology,
called the \textbf{strong topology}, whose basis is the family of  open sets of the \emph{gts}.
(The members of the strong topology are exactly weakly open subsets of the space.)
Passing to the strong topology forms a functor $()_{top}:\lss \to \mathbf{Top}$.

A locally small space is called: \textbf{paracompact} if there is a locally finite covering of the space by small open subsets, and  \textbf{Lindel\"of} if there is a countable admissible covering of the space by small open subsets. Each connected paracompact locally small space is 
Lindel\"of (the proof of I.4.17 in \cite{LSS} is purely topological).

\begin{prop}[cf. I.4.6 of \cite{LSS}]
For each  paracompact locally small space $X$,
 the weak closure $\overline{Y}$ (that is: the closure in the strong topology) of a small set $Y$ is small. 
\end{prop}
\begin{proof}
Take a locally finite covering $\mc{U}$ of $X$ by small open subsets.
The set $Y$ is covered by a finite subcover $\mc{U}_0$ of $\mc{U}$ of all members of $\mc{U}$
that meet $Y$.
Then $\overline{Y}$ and the union of $\mc{U}\setminus \mc{U}_0$ are disjoint,
 and $\overline{Y}$ is contained in the union of $\mc{U}_0$, which is a small set. 
\end{proof}

\begin{exam}
If a metric space $(X,d)$ satisfies the ball property

(\textbf{BP}) each intersection of two open balls is a finite union of open balls,

then $X$ has  a natural generalized topology, where an open set $Y$ is such a subset of $X$ that
the trace of $Y$ on each open ball is a finite union of open balls, and the admissible coverings
are such open coverings that are essentially finite on each open ball.
Then open balls are small sets. The covering of the space by all open balls is admissible (each small set is covered by one open ball), and $X$ is a locally small space satisfying the first axiom of countability. The family of all open balls is a basis of the generalized topology.
\end{exam}

A subset $Y$ of a locally small space $X$ is \textbf{locally constructible} if 
each intersection $Y\cap U$ with a small open $U\subseteq X$ is constructible in $U$ (so also in $X$).
The Boolean algebra of locally constructible subsets of a locally small space may be strictly larger than the Boolean algebra of constructible subsets (to see this one can construct a sequence $X_n$ of constructible subsets of some small spaces $Z_n$ each $X_n$ needing
at least   $n$ open sets in the description, and then glue the spaces $Z_n$ into one locally small space).

We will say that a locally small space has the \textbf{closure property} if the following holds:

(\textbf{CP}) the weak closure of a small locally closed  subset is a closed subset.
\begin{f}
Each locally small space with a regular strong topology and with the closure property is weakly regular.
\end{f}
\begin{f}\label{48}
Each locally constructible subset of a locally small space is a locally essentially finite
union of small locally closed subsets. 
\end{f}

\begin{prop} If a locally small space has the closure property, then the weak closure of each locally constructible set is a closed set.
\end{prop}
\begin{proof}
By the Fact \ref{48}, the weak closure of each locally constructible set is a
union of weak closures of some locally essentially finite family of small constructible sets.
But the family of these closures is also locally essentially finite, and the thesis follows from Proposition \ref{closed}. 
\end{proof}

By the above, if a locally small space has the closure property, then the closure operator of the generated topology restricted to the class of locally constructible sets may be considered as the \textit{closure operator} of the generalized topology (in general, the closure operator on a \emph{gts} does not exist).

Notice that each subset $Y$ of a locally small space $X$ forms a (locally small) subspace, since if $X=\stackrel{a}{\bigcup}_{\alpha} X_{\alpha}$ with all $X_{\alpha}$ small open, then the formula $Y=\stackrel{a}{\bigcup}_{\alpha} (X_{\alpha}\cap Y)$ defines a locally small space that is a subspace of $X$.
Hence Proposition \ref{cu}, the concept of a connected component, and a quasi-component extend to the category $\lss$.
Also the equalizers for parallel pairs of morphisms exist in $\lss$.
Together with Lemma \ref{lokmal}   this gives

\begin{thm} The category $\lss$ is finitely complete.
\end{thm}

\section{Weakly small spaces}
Now we start to reintroduce the theory of weakly semialgebraic spaces from
\cite{WSS} on a topological level.

A \textbf{weakly} (or \textbf{piecewise}) \textbf{small space} is a \emph{gts} $X$ having a family $(X_{\alpha})_{\alpha \in A}$
of  closed small subspaces indexed by a partially ordered set $A$ such that the following conditions hold:

W1) $X$ is the union of all $X_{\alpha}$'s as sets,

W2) if $\alpha \leq \beta$ then $X_{\alpha}$ is a (closed, small) subspace of $X_{\beta}$, 

W3) for each $\alpha \in A$ there are only finitely many $\beta \in A$ such that $\beta  < \alpha$,

W4) for each two $\alpha,\beta \in A$ there is $\gamma \in A$ such that $X_{\alpha}
\cap X_{\beta} = X_{\gamma}$,

W5) for each two $\alpha,\beta \in A$ there is $\gamma \in A$ such that $\gamma \geq \alpha$ and $\gamma \geq \beta$,

W6) the \emph{gts} $X$ is the inductive limit of the directed family $(X_{\alpha})_{\alpha \in A}$, which means:

a) a subset $U$ of $X$ is open iff all sets $X_{\alpha}\cap U$ are open in respective $X_{\alpha}$s,

b) an open family $\mc{U}$ is admissible iff all for each $\alpha \in A$ the family $\mc{U}\underline{\cap} X_{\alpha}$ is admissible 
(=essentially finite) in respective $ X_{\alpha}$.

Such a family $(X_{\alpha})_{\alpha \in A}$ is called an \textbf{exhaustion} of $X$.
The weakly small spaces form a full subcategory $\wss$ of $\gts$.
Here ``admissible'' means ``piecewise essentially finite'', where ``piecewise'' means ``when restricted to a member of the exhaustion (chosen to witness that $X$ is a weakly small space)''.
Members of this exhaustion of $X$ may be called \textbf{pieces}.
If $(X_{\alpha})_{\alpha \in A}$ is  an exhaustion of $X$, then we will write
$X=\stackrel{e}{\bigcup}_{\alpha \in A} X_{\alpha}$.
 
The \textbf{index function} for the exhaustion  $(X_{\alpha})_{\alpha \in A}$ is
the function $\eta: X\rightarrow A$ given by the formula
$$ \eta(x)=\inf \{ \alpha \in A \mid \: x\in X_{\alpha} \}.$$
Here infimum exists thanks to W3).
The index function $\eta$ gives a decomposition of the space $X$ into small locally closed subspaces $X^0_{\alpha}=\eta^{-1}(\alpha)=X_{\alpha}\setminus \bigcup_{\beta <\alpha} X_{\beta}$.

\begin{f}
A subset $Y$ of a weakly small space $X$ is closed if and only if it is piecewise closed.
\end{f}
\begin{prop}
A piecewise essentially finite union of closed subsets of a weakly small space is closed.
\end{prop}
\begin{proof}
An essentially finite union of closed subsets is closed. Hence a piecewise essentially finite union of closed sets is piecewise closed, thus closed.
\end{proof}

Notice that the chosen exhaustion of a weakly small space is a piecewise essentially finite (relative to this exhaustion) family of closed sets, and remind that
a constructible subset of a piece is a finite union of locally closed subsets of a piece.
A \textbf{weakly} (or \textbf{piecewise}) \textbf{constructible subset} is such a subset $Y\subseteq X$
that has constructible intersections with all members of the chosen exhaustion $(X_{\alpha} )_{\alpha \in A}$.

\begin{prop}
Piecewise constructible subsets of a weakly small space are exactly piecewise essentially finite unions of locally closed subsets of pieces.
\end{prop}
\begin{proof}
Piecewise constructible subsets are piecewise essentially finite unions of locally closed 
subsets of pieces. An essentially finite union of locally closed subsets of pieces
is a finite union of locally closed subsets of a single piece, so a constructible subset of a piece. Now apply ``piecewise''.
\end{proof}

The \textbf{strong topology} on $X=\stackrel{e}{\bigcup}_{\alpha \in A}X_{\alpha}$ is the topology that makes the topological space $X$ the respective inductive limit of the system of topological spaces $X_{\alpha}$. Its members are all the piecewise weakly open subsets, not only the weakly (piecewise) open subsets. Hence the open sets from the generalized topology may not form  a basis of the strong topology (see Appendix C of \cite{WSS}).
Another unpleasant fact about the weakly small spaces (comparising with the locally small spaces)  is that points may not have small neighborhoods (consider an infinite wedge of circles as in Example 4.1.8 of \cite{WSS}).

We will say that a weakly small space has the \textbf{closure property} if the following holds:

(\textbf{CP}) the weak closure  of a locally closed subset of a piece is a closed subset.

Notice that the strong topology and the generated topology coincide on every piece.
If the space $X$ has the closure property, then
the topological closure operator restricted to the class of  constructible subsets of pieces
of $X$ may be treated as the \textit{closure operator} of the generalized topology.

Notice that each subset $Y$ of a weakly small space $X$ forms a (weakly small) subspace, since if $X=\stackrel{e}{\bigcup}_{\alpha} X_{\alpha}$ with all $X_{\alpha}$ small closed, then $Y=\stackrel{e}{\bigcup}_{\alpha} (X_{\alpha}\cap Y)$ defines a weakly small space that is a subspace of $X$.
Hence Proposition \ref{cu}, the concept of a connected component, and a quasi-component extend to the category $\wss$.

\begin{thm}[cf. IV.2.1 in \cite{WSS}]\label{scp}
If a weakly small space  $X=\stackrel{e}{\bigcup}_{\alpha \in A} X_{\alpha}$ is strongly $T_1$, and $L$ is a small space, then
 for each  strictly continuous mapping $f:L\rightarrow X$ there is $\alpha_0 \in A$ such that $f(L)\subseteq X_{\alpha_0}$.
\end{thm}
\begin{proof}
Let $\eta:X\rightarrow A$ denote the index function for the exhaustion $(X_{\alpha})_{\alpha \in A}$.
If $\eta(f(L))$ were infinite, then for each $\alpha \in \eta(f(L))$ we could choose $x_{\alpha}\in f(L)$ with $\eta(x)=\alpha$ and some $y_{\alpha}\in f^{-1}(x_{\alpha})$.
Set $S=\{x_{\alpha}: \alpha \in \eta(f(L))\}$.
For each $\gamma \in A$, the set $S\cap X_{\gamma}$ is finite. Since $X$ is strongly $T_1$,
the set $S$ is closed as well as each of its subsets, so $S$ is topological discrete. Then $\{ y_{\alpha}:  \alpha \in \eta(f(L))  \}\subseteq f^{-1}(S)$ is a topological discrete, infinite, and small subset of $L$. This is a contradiction.

Hence $\eta(f(L))$ is finite, and there is $\alpha_0\geq \eta(f(L))$. 
We get $f(L)\subseteq X_{\alpha_0}$.
\end{proof}

\begin{thm}[cf. IV.2.2 in \cite{WSS}]
If a weakly small space $X$ with an exhaustion $(X_{\alpha})_{\alpha \in A}$ is strongly $T_1$, then each small subspace $L$ of $X$ is contained in some $X_{\alpha_0}$. In particular, each member  $X_{\beta}$ of any exhaustion  $( X_{\beta})_{\beta \in B}$ of $X$ is contained in some member $X_{\alpha_0}$ of the initial exhaustion.
\end{thm}
\begin{proof}
For  each such $L$, the inclusion mapping $i:L\rightarrow X$ is strictly continuous. By Theorem \ref{scp}, the set $i(L)=L$ is contained in a member of the exhaustion $(X_{\alpha})_{\alpha \in A}$.
\end{proof}

Let us denote by $\wss_1$ the full subcategory of $\wss$ composed of  strongly $T_1$ objects of $\wss$. In this category the term ``piecewise'' does not depend on an exhaustion (as it may be expressed by ``when restricted to a closed small set''), hence 
passing to the strong topology forms a functor $()_{stop}: \wss_1\to \mathbf{Top}$.

\begin{thm}\label{wss1}
The category $\wss_1$ is finitely complete.
\end{thm}
\begin{proof} It is enough to consider binary products.
If $X$, $Y$ have exhaustions $(X_{\alpha})_{\alpha \in A}$, $(Y_{\beta})_{\beta \in B}$, respectively, then $( X_{\alpha}\times Y_{\beta})_{(\alpha,\beta)\in A\times B}$ is an exhaustion defining the weakly small space $X\times Y$.  Notice that $X\times Y$ is
strongly $T_1$.  The projections are obviously strictly continuous, and for 
strictly continuous $f:Z\to X$, $g:Z\to Y$, the mapping $(f,g):Z\to X\times Y$ is strictly continuous, since if $Z=\stackrel{e}{\bigcup}_{\gamma\in \Gamma} Z_{\gamma}$, then for each $Z_{\gamma_0}$ there are $X_{\alpha_0}$ and $Y_{\beta_0}$ such that $f(Z_{\gamma_0})\subseteq X_{\alpha_0}$ and $g(Z_{\gamma_0})\subseteq Y_{\beta_0}$. Thus it is the product of $X$ and $Y$ in $\wss_1$. 

Since each subset forms a subspace, the existence of equalizers for pairs of
parallel mappings is clear.
\end{proof}

\begin{f}
If $X$ is a small space and $Y$ is an object of $\wss_1$, then the projection $X\times Y\rightarrow Y$ is an open and closed mapping.
\end{f}

\section{Spaces over structures}
In this chapter, we deal with locally definable and weakly definable spaces over mathematical structures.

Assume that $\mc{M}$ is any (one sorted, first order for simplicity) structure in the sense of model theory.
A \textbf{function sheaf over $\mc{M}$} on a \emph{gts}
$X$ is a sheaf $F$ of sets on $X$ (the sheaf property is assumed only for admissible coverings) such that for each open $U$  
the set $F(U)$ is contained in the set $M^U$ of all functions
 from $U$ into $M$, and the restrictions of the sheaf
 are the set-theoretical restrictions of functions.
A \textbf{space over $\mc{M}$} is a pair $(X, O_X )$, where $X$ is a
 \emph{gts}  and $O_X$ is  a function sheaf over $\mc{M}$ on $X$.
A \textbf{morphism} $f:(X, O_X )\rightarrow (Y, O_Y)$ of spaces over $\mc{M}$ is a  strictly continuous mapping $f: X\rightarrow Y$ such that for each open subset $V$ of $Y$
the  set-theoretical substitution $h \mapsto h\circ f$ gives the mapping $f^{\#}_V:O_Y (V) \rightarrow O_X (f^{-1} (V))$. (We could informally say that $f^{\#} : O_Y \rightarrow O_X $ is the ``morphism of function sheaves'' over $\mc{M}$ induced set-theoretically by $f$.
We can also define for function sheaves
$$(f_{*}O_X) (V)=\{ h:V\to R |\: h\circ f\in O_X(f^{-1}(V))\},$$
and then each $f^{\#}_V:O_Y(V)\to f_{*}O_X (V)$ is an inclusion.)
An \textbf{isomorphism} is an invertible morphism. We get a category $\spa(\mc{M})$ of spaces over $\mc{M}$ and their morphisms.

\begin{thm}[cf. I.2, p. 11, in \cite{LSS}]
For each structure $\mc{M}$, the  category $\spa(\mc{M})$ is cocomplete.
\end{thm}
\begin{proof}
Assume that $X$ is the inductive limit  (in the category $\gts$) of a diagram $D$ with objects $(X_{i})_{i\in I}$, indexed by a small category $I$, and the canonical morphisms $\phi_i: X_i\rightarrow X$.

Assume additionally that function sheaves $O_{X_i}$ over $\mc{M}$   are given.
Then define
$$ O_X (U)=\{ h:U\rightarrow M \mid \: \forall i\in I \:\: h\circ \phi_i\in O_{X_i} (\phi_i^{-1}(U))\}.$$
Since all $O_{X_i}$ are function sheaves and the section amalgamation in function sheaves is given by the union of graphs of sections, also $O_X$ is a function sheaf. By the above definition, each  $\phi_i$ is a morphism of spaces over $\mc{M}$.

If another space $(Y,O_Y)$ with a set of morphisms $\psi_i:X_i\rightarrow Y$ is given, 
then there is a unique morphism $\eta:X\rightarrow Y$ in $\gts$ such that $\psi_i = \eta\circ\phi_i$ for all $i$. But $\eta$ is also a morphism in $\spa(\mc{M})$, because
the condition $k\circ \eta\in O_X(\eta^{-1}(W)) $
for each $k\in O_Y(W)$ is visibly satisfied, since, for each $i\in I$, we have
$$k\circ\eta\circ\phi_i=k\circ\psi_i\in O_{X_i}(\psi_i^{-1}(W))=O_{X_i}(\phi_i^{-1}(\eta^{-1}(W)))     .$$  
\end{proof}


Now assume that a topology is given on the underlying set $M$ of a (first order) structure $\mc{M}$. We explicitely demand that the product topologies on cartesian powers $M^n$ should be considered.
We will call such $\mc{M}$ a \textbf{weakly topological structure}. (This setting seems to coincide with the case (i) in the introduction of \cite{Pi}. We do not
explore a special language $L_t$ for topological structures considered in the case (ii) of the introduction of \cite{Pi} or in \cite{FZ}.)
Then all projections $\pi_{n,i}:M^n\to  M^{n-1}$ (forgetting the $i$-th coordinate) are continuous and open. Thus, for example, the field of complex numbers $(\mb{C},+,\cdot)$
considered with the euclidean topology (but not with the Zariski topology) is a weakly topological structure. Also the fields $\mb{Q}_p$ of p-adic numbers considered with their natural topologies (coming from valuations) are weakly topological. 

Now  for each definable set $D\subseteq M^n$, we  set (as in \cite{ap}): 

a) an \emph{open subset} of a \emph{gts} means a relatively open, definable subset;

b) an \emph{admissible family} of a \emph{gts} means an essentially finite family.

Each such $D$ becomes a small \emph{gts}, and $\pi_{n,i}$ become strictly continuous and open.
We define the structure sheaf $O_D$ as the sheaf of all definable continuous functions
from respective (\emph{gts}-)open subsets $U\subseteq D$  into $M$.
Thus $(D, O_D)$ becomes a space over $\mc{M}$.

An \textbf{open subspace} of a space over $\mc{M}$ is an open subset
 of its \emph{gts} together with the function sheaf of the space
  restricted to this open set.
A \textbf{small subspace} $K$ of a space $X$ over $\mc{M}$ is a small subset of its \emph{gts} with the function sheaf $O_K$ composed of all finite open unions of restrictions of sections of the function sheaf $O_X$  to the relatively open subsets of the small subspace $K$.

 An \textbf{affine definable space} over $\mc{M}$ is a space over $\mc{M}$
isomorphic to a definable subset of some $M^n$ considered with its usual structure of a space over $\mc{M}$.
A \textbf{definable space} over $\mc{M}$ is a space over $\mc{M}$ that has
 a finite open covering  by affine definable subspaces. (The structure sheaf is determined  in an obvious way.)
A \textbf{locally definable space} 
 over $\mc{M}$ is a space over $\mc{M}$ 
that has an admissible covering  by affine definable open subspaces.
(The sections of a structure sheaf are admissible unions of sections of the structure sheaves
of affine open subspaces.)
A \textbf{weakly definable space} $X$ over $\mc{M}$ is a space over $\mc{M}$ that has an exhaustion $(X_{\alpha})_{\alpha\in A}$
composed of  definable (small) subspaces such that a function $h:V\rightarrow M$ is a section of $O_X$ iff
all restrictions $h|_{V\cap X_{\alpha}}$ are sections of respective $O_{X_{\alpha}}$s. 

\textbf{Morphisms} of affine definable spaces, definable spaces,  locally
 definable spaces, and weakly definable spaces over $\mc{M}$ are their morphisms as spaces over $\mc{M}$. We get
 full subcategories $\ads(\mc{M})$, $\ds(\mc{M})$, $\lds(\mc{M})$, $\wds(\mc{M})$ of
 $\spa(\mc{M})$.

Notice that each definable subset of an affine definable space is also an affine definable space. Thus definable subsets of definable spaces are naturally \textbf{subspaces} in $\ds(\mc{M})$. Locally definable subsets of locally definable spaces are naturally  \textbf{subspaces} in $\lds(\mc{M})$. Also piecewise definable subsets of weakly definable spaces are naturally  \textbf{subspaces} in $\wds(\mc{M})$. It is clear that 
$\ads(\mc{M})$, $\ds(\mc{M})$, $\lds(\mc{M})$, $\wds(\mc{M})$ have equalizers for pairs of parallel morphisms.

The notions of \textbf{paracompactness} and \textbf{Lindel\"ofness} for locally definable spaces coincide with
their counterparts for the underlying locally small spaces.
 
As in \cite{ap}, we get from the definitions

\begin{fs}[Facts 5 in \cite{ap}]  Definable spaces are small. Every small  subspace of a locally definable space is definable. 
\end{fs}

Notice that if $M$ is $T_1$, then all objects of $\lds(\mc{M})$ and $\wds(\mc{M})$ are strongly $T_1$.  Moreover, each weakly $T_1$ locally definable or weakly definable space over $\mc{M}$ is strongly $T_1$. (It is visible for affine spaces, and extends to the general case by 
applying ``locally'' or ``piecewise''.) Thus we can speak just about $T_1$ objects of
 $\lds(\mc{M})$ or of $\wds(\mc{M})$ for any $\mc{M}$. 
 From Theorem \ref{scp}, we get
\begin{f}
If a weakly definable space is  $T_1$, then each small subspace 
is contained in a piece, hence it is a definable subspace.
\end{f}

\begin{thm}
For each weakly topological structure $\mc{M}$, the  products of small spaces in the category $\spa(\mc{M})$ exist.
\end{thm}
\begin{proof}
Let a family $\{(X_{\alpha}, O_{X_{\alpha}})\}_{\alpha \in A}$ of small spaces over $\mc{M}$ be given. By the Proposition \ref{gtsprod}, the
family $\{ X_{\alpha}\}_{\alpha \in A}$ has the product $X$ in $\gts$.
Now, for any $h\in O_{X_{\alpha}}(V)$, the function $h\circ \pi_{\alpha}:\pi_{\alpha}^{-1}(V)\rightarrow M$ may be called a ``function of one variable'', where $\pi_{\alpha}:X\rightarrow X_{\alpha}$ is the projection.
Define sections of $O_X$ as all admissible (so essentially finite) unions of compatible functions of one variable.
Then the projections $\pi_{\alpha}$ are morphisms of $\spa(\mc{M})$.
For a family $f_{\alpha}: Y\rightarrow X_{\alpha}$ of morphisms of $\spa(\mc{M})$,
the induced mapping $(f_{\alpha})_{\alpha \in A}:Y\rightarrow X$ is a morphism of
 $\spa(\mc{M})$, since all the sections of $O_X$ are finite unions of the functions of one variable, which are  determined by sections of $O_{X_{\alpha}}$s.
\end{proof}

Unfortunately, even finite products of affine definable spaces in
$\spa(\mc{M})$ usually are not affine definable spaces.

\begin{prop}
For each weakly topological structure $\mc{M}$, the categories $\ads(\mc{M})$ and $\ds(\mc{M})$ are finitely complete.
\end{prop}
\begin{proof}
We need only to check if finite products exist in $\ads(\mc{M})$ and 
$\ds(\mc{M})$. For definable sets $D\subseteq M^n$ and $E\subseteq M^m$, the product is the definable set
$D\times E\subseteq M^{n+m}$ (with its natural space structure). Indeed, the projections are continuous definable mappings, and for any continuous definable $f:Z\to D$, $g:Z\to E$, with $Z$ definable in some $M^k$, the induced mapping $(f,g):Z\to D\times E$ is continuous definable, since the topology of $D\times E$ is the product topology.
For a continuous definable function $h\in O_{D\times E}(U)$, the function $h\circ (f,g):(f,g)^{-1}(U)\to M$ is continuous definable.
This finishes the proof for the affine definable case. The definable case is similar.
\end{proof}
Notice that the canonical projections from a finite product of affine definable spaces to its factors are open morphisms, but not in general closed morphisms. 

\begin{f}[cf. I.1.3 in \cite{LSS}]
A mapping $f:X\to Y$ between objects of $\lds(\mc{M})$ is a morphism iff the image 
of each open definable subspace  is contained in an open definable subspace  and when restricted to open definable subspaces in the domain and in the range, the mapping is continuous definable.
\end{f}


\begin{thm}[cf. I.2.5 in \cite{LSS}]\label{fplds}
For each weakly topological structure $\mc{M}$, the category $\lds(\mc{M})$
is finitely complete.
\end{thm}
\begin{proof}
We need only to check the existence of finite products.
For locally definable spaces $(X,O_X)=\stackrel{a}{\bigcup}_{\alpha\in A} (X_{\alpha},O_{X_{\alpha}})$, $(Y,O_Y)=\stackrel{a}{\bigcup}_{\beta\in B} (Y_{\beta},O_{Y_{\beta}})$, where 
 $(X_{\alpha},O_{X_{\alpha}})$, $(Y_{\beta},O_{Y_{\beta}})$ are affine definable subspaces over $\mc{M}$, consider the space $$\stackrel{a}{\bigcup}_{(\alpha,\beta)\in A\times B} 
(X_{\alpha},O_{X_{\alpha}})\times (Y_{\beta},O_{Y_{\beta}})=\stackrel{a}{\bigcup}_{(\alpha,\beta)\in A\times B} (X_{\alpha}\times Y_{\beta}, O_{X_{\alpha}\times Y_{\beta}}).$$
This is the product of given locally definable spaces, since the projections are visibly morphisms in
$\lds(\mc{M})$, and for  given morphisms $f:Z\to X$, $g:Z\to Y$ in $\lds(\mc{M})$ from the space
$(Z,O_{Z})=\stackrel{a}{\bigcup}_{\gamma\in \Gamma} (Z_{\gamma},O_{Z_{\gamma}})$,
the induced mapping $(f,g):Z\to X\times Y$ is a morphism in  $\lds(\mc{M})$. Indeed, we may assume $\{Z_{\gamma}\}_{\gamma \in \Gamma}$ is
an admissible covering by affine definable spaces that is a refinement of the preimage of the covering
$\{X_{\alpha}\times Y_{\beta}\}_{(\alpha,\beta)\in A\times B}$ and thus each $Z_{\gamma}$ is mapped into some $X_{\alpha}\times Y_{\beta}$ by some continuous definable mapping, and the whole $(f,g)$ is an admissible union of such partial mappings. 
\end{proof}

We can also generalize the following concept from \cite{ap}: a subset $Y$ of a locally definable space $X$ (over $\mc{M}$) may be called \textbf{local} if for each $y\in Y$ there is an open small neighborhood  $U$ of $y$ such that $Y\cap U$ is definable in $U$. 
On such $Y$, we could define a locally definable space by the formula
$$Y=\stackrel{a}{\bigcup} \{ Y\cap U\mid U\mbox{ small, open}, Y\cap U\mbox{ definable in }U \} .$$
(This definition is canonical in the sense that it does not depend on any arbitrary choice of open neighborhoods.) However, the use of local subsets does not reflect this
locally definable space (see Examples 11 and 12 in \cite{ap}). Any weakly open and any weakly discrete set is local.

Let us denote by $\wds_1(\mc{M})$ the full subcategory of $T_1$ objects of $\wds(\mc{M})$ for any $\mc{M}$.

\begin{f}[cf. IV.2.3 in \cite{WSS}]
A mapping $f:X\to Y$ between  objects of $\wds_1(\mc{M})$ is a morphism iff the image of a piece is contained a in piece and when restricted to pieces in the domain and in the range, the mapping is continuous definable.
\end{f}

Notice that for two morphisms $f_1:X_1\to Y_1$, $f_2:X_2\to Y_2$ of $\lds(\mc{M})$ or of $\wds_1(\mc{M})$, their cartesian product $f_1\times f_2:X_1\times X_2\to Y_1\times Y_2$ is a morphism. 

\begin{thm}[cf. IV.3, p. 32, in \cite{WSS}]
For each weakly topological structure $\mc{M}$, the category $\wds_1(\mc{M})$ is finitely complete.
\end{thm}
\begin{proof}
We need only to check the existence of finite products.
Notice that  by Theorem \ref{fplds}, finite products exist in the category $\ds(\mc{M})$, and the product of $T_1$ spaces is $T_1$.
 We will drop the structure sheaves in notation.
If $(X_{\alpha})_{\alpha \in A}$, $(Y_{\beta})_{\beta \in B}$ are exhaustions of objects $X$, $Y$, then define 
$X\times Y= \stackrel{e}{\bigcup}_{(\alpha,\beta)\in A\times B} ( X_{\alpha}\times Y_{\beta})_{(\alpha,\beta)\in A\times B}$. This is  an exhaustion determining the structure sheaf.  The resulting space is clearly $T_1$, and the projections $X\times Y\to X$, $X\times Y\to Y$ are visibly morphisms of 
$\wds_1(\mc{M})$. If $Z=\stackrel{e}{\bigcup}_{\gamma\in \Gamma} Z_{\gamma}$ is any object of 
$\wds_1(\mc{M})$ and $f:Z\to X$, $g:Z\to Y$ are morphisms,then  the mapping $(f,g):Z\to X\times Y$ is a morphism because:  we may assume that for each 
 $Z_{\gamma_0}$ there are $X_{\alpha_0}$ and $Y_{\beta_0}$ such that $f(Z_{\gamma_0})\subseteq X_{\alpha_0}$ and $g(Z_{\gamma_0})\subseteq Y_{\beta_0}$, the generalized topology of the product of definable spaces contains the generalized topology of their (generalized topological) product, and the definable continuous mappings are stable under composition and juxtaposition.
 Thus $X\times Y$ is the product of $X$ and $Y$ in $\wds_1(\mc{M})$.
\end{proof}

\begin{prop}[cf. I.3.5 in \cite{LSS}]
For any weakly topological structure $\mc{M}$, the category $\lds(\mc{M})$
has fiber products.
\end{prop}
\begin{proof}
The fiber product of morphisms $f:X\to Z$ and $g:Y\to Z$ is 
the subspace $(f\times g)^{-1}(\Delta_Z)\subseteq X\times Y$, where $\Delta_Z$ is the diagonal of the space $Z$ (locally definable in $Z\times Z$).
\end{proof}
 
\begin{prop}[cf. IV.3.20 in \cite{WSS}]
For any weakly topological structure $\mc{M}$, the category $\wds_1(\mc{M})$
has fiber products.
\end{prop}
The proof of this proposition is similar to the proof of the previous proposition.

Let $\mc{C}$ be  one of the categories: $\lds(\mc{M})$ or $\wds_1(\mc{M})$. 
We will say that an object $Z$ of $\mc{C}$ is $\mc{C}$\textbf{-complete} if the mapping 
$Z\to \{*\}$ is universally closed in $\mc{C}$, which means that for each object $Y$ of $\mc{C}$ the projection $Z\times Y\to Y$ (which is the base extension of $Z\to \{*\}$) is a closed mapping. 

\begin{prop}\label{clcompl}
A closed subspace of a $\mc{C}$-complete space is $\mc{C}$-complete.
\end{prop}
\begin{proof}
Let $C$ be a closed subspace of a $\mc{C}$-complete space $Z$. For any $Y$ in $\mc{C}$, the space $C\times Y$ is a closed subset of $Z\times Y$. Thus the image under the projection along $C$ of any closed subset of $C\times Y$ is closed in $Y$.
\end{proof}

\begin{prop}\label{clsub}
Every weakly closed subspace of an object of $\mc{C}$  is a closed subspace.
\end{prop}
\begin{proof}
This is clear for affine definable spaces.
In general, closedness may be checked ``locally'' or ``piecewise''.
\end{proof}

\begin{prop}\label{compl}
The image $g(C)$ of a $\mc{C}$-complete subspace $C$ of $Y$
 under a morphism $g:Y\to Z$ in $\mc{C}$ is a $\mc{C}$-complete subspace.
\end{prop}
\begin{proof}
Consider an object $W$ of $\mc{C}$. For a closed subset $A$ of $g(C)\times W$, we have $\pi_{g(C)}(A)=\pi_C ((g|_C\times id_W)^{-1}(A))$ is a closed set.
\end{proof}

\begin{prop}\label{hausdorff}
For an object $Z$ of $\mc{C}$, the following conditions are equivalent:

a)  $Z$ is weakly Hausdorff;

b) $Z$ is strongly Hausdorff;

c) $Z$  has its diagonal $\Delta_Z$ closed.
\end{prop}
\begin{proof}
If a space  $X$ is weakly Hausdorff, then its generated topology is Hausdorff, thus the
diagonal $\Delta_X$ is closed in the generated topology. But $\Delta_X$ is always a subspace, thus, by Proposition \ref{clsub}, a closed subspace.

If the diagonal  $\Delta_X$ of a space $X$ is a closed subspace, then the generated topology of $X$ is Hausdorff, so $X$ is weakly Hausdorff. But it is also strongly $T_1$, since for any $x_0\in X$, the set $\{(x_0,x_0)\}$ is closed in $\{ x_0\}\times X$ and the projection 
$\{ x_0\}\times X \to X$ is closed. Thus $X$ is strongly Hausdorff.
\end{proof}

Notice that if the topology of $M$ is Hausdorff, then:

a) all $M^n$ as well as all objects of $\ads(\mc{M})$ are (strongly) Hausdorff;

b)  all objects of $\wds(\mc{M})$ having exhaustions consisting of affine definable spaces are (strongly) Hausdorff;

c) an object of  $\ds(\mc{M})$ may not be (weakly) Hausdorff.

Because of \ref{hausdorff}, we may speak just about Hausdorff objects of $\mc{C}$.

\begin{prop} \label{complete}
If $Y$ is a Hausdorff  object  of $\mc{C}$, then:

a) each of its $\mc{C}$-complete subspaces  is closed;

b) the graph of each morphism $f:X\to Y$ in $\mc{C}$ is closed;
\end{prop}
\begin{proof}
a) If $C$ is a $\mc{C}$-complete subspace of $Y$, then $\Delta_Y\cap (C\times Y)$ is relatively closed in $C\times Y$, and its projection on $Y$, equal to $C$, is closed.

b) The graph of $f$ is a subspace of $X\times Y$ being the preimage of $\Delta_Y$
by a morphism $f\times id_Y$. 
\end{proof}

Suppose $\mc{M}$ satisfies the condition\\
$(\mathbf{DCCD})$ $\mc{M}$ admits (finite) cell decomposition with definably connected cells\\
(here we assume only that a cell is a definable set, for stronger notions of a cell
see \cite{Ma} or \cite{Sch}), then:

a) each object of $\ds(\mc{M})$ has a finite number of clopen definable connected components, and is a finite (generalized) topological direct sum of them;
 
b) each object of $\lds(\mc{M})$ is a (locally finite) direct sum
of its connected components;
 
c) each object of $\wds(\mc{M})$ is a piecewise finite union of its clopen connected components, so also a direct sum of its connected components.

Thus If $\mc{M}$ satisfies $(DCCD)$, then all objects of $\lds(\mc{M})$ and of $\wds(\mc{M})$ satisfy $(ACC)$. 
All o-minimal structures $\mc{M}$ (with their  order dense) satisfy $(DCCD)$, but t-minimality in the sense of \cite{Sch} does not guarantee $(DCCD)$. 

\begin{lem} \label{disncomp}
If $\mc{C}$ has a $T_1$ object having a non-closed countable set, then
a topological discrete infinite space is not $\mc{C}$-complete.
\end{lem}
\begin{proof}
Let $Z$ be an infinite discrete space, and $W$ a $T_1$ object with a non-closed countable set $C=\{ c_n \}_{n\in \mb{N}}$
(given by an injective sequence). Take an injective sequence
$\{ z_n \}_{n\in \mb{N}}$ of elements of $Z$.
Then $C=\pi_Z (\{ (z_n,c_n):n\in \mb{N}    \})$, so $Z$ is not complete. 
\end{proof}

\begin{thm}[cf. I.5.8 in \cite{LSS}]\label{fncc}
If  
there is a $T_1$ object $Z$ of $\mc{C}$ having a non-closed countable subset, then each  $\mc{C}$-complete space 
intersects only a finite number of components of a generalized topological sum.
\end{thm}

\begin{proof}
The $\mc{C}$-complete space $C$ is a direct sum of the family of its connected components. Take one point from each connected component. The resulting subspace $S$ is discrete and closed, so complete by Proposition \ref{clcompl}. If it were infinite, then we could assume (by taking a closed subspace if necessary) it is  countable.   An injection from $S$ to  a countable non-closed subset of some Hausdorff object of $\mc{C}$ would give  a counterexample of  c) in Proposition \ref{complete}. Hence $S$ is finite.
\end{proof}

\begin{thm}[cf. I.5.10 in \cite{LSS}]
Assume that 
there is a $T_1$ object $Z$ of $\mc{C}$ having a non-closed countable subset.
Let $C$ be a $\mc{C}$-complete space. 
Additionally, if $\mc{C}=\lds(\mc{M})$ then $C$ is assumed to be $T_1$
and Lindel\"of.
Then $C$ is a definable space.
\end{thm}
\begin{proof}
In the case of a locally definable space: 
If $C=\stackrel{a}{\bigcup}_{n\in \mb{N}} C_n$ is an admisible covering by (affine) definable spaces and $C$ is not definable, then we may assume that for each $n\in \mb{N}$ we can choose 
$x_n\in C_{n}\setminus (C_{0}\cup ...\cup C_{{n-1}})\neq \emptyset$. The set $B=\{ x_n |n\in \mb{N}\}$ is infinite, not small, locally finite subspace, so, as a $T_1$ space, topological discrete.
Moreover, $B$ is closed, thus complete (by Proposition \ref{clcompl}).

In the case of a weakly definable space: if $C$ is not definable, then  the index function $\eta$ of an exhaustion
$(C_{\alpha})_{\alpha\in A}$ has infinite image, and
we can choose an element $x_{\alpha}\in C^0_{\alpha}$ for each $\alpha\in \eta(C)$.
The set $B=\{ x_{\alpha}: \alpha\in \eta(C) \}$ is an infinite (but piecewise finite) closed subspace all of whose subsets are also closed. Thus $B$ is a complete (by Proposition \ref{clcompl}) and topological discrete space.
We may assume $B$ is countable.

We get a contradiction with Lemma \ref{disncomp}. 

\end{proof}

Remind (after \cite{Pi}) that a structure $\mc{M}$ is a  \textbf{first order topological structure} (called in \cite{FZ} a topological structure with explicitely definable topology) if the basis of the topology on $M$ is uniformly definable in $\mc{M}$. (There is a formula $\Phi(x,\bar{y}) $ of the (first order) language of $\mc{M}$ such that the family
$\{ \Phi(x,\bar{a})^M : \bar{a}\subseteq M\}$ is the basis of the topology of $M$.)   
If $\mc{M}$ is a first order  topological structure, then the relative closure 
of a definable subset in a definable set is definable (see \cite{Pi}). 
Thus the closure operator for the generalized topology exists for the locally definable
subsets of locally definable spaces, and for definable subsets of pieces in weakly definable spaces over $\mc{M}$.

\section{Open problems}
Generalized topology in the sense of H. Delfs and M. Knebusch is a new chapter in general topology. That is why there are many unanswered questions
in this topic. The following questions are suggested to the reader:

1) Does any subset of a \emph{gts} form a subspace?  What kinds of subsets
(other than the open subsets and the small subsets) always form subspaces?
For what classes of generalized topological spaces (other than locally small or weakly small) all subsets form subspaces?

2) Does \textbf{GTS} have products?

3)  For a family $\{ V_{i}\}_{i\in I}$ of subsets, and a family $\{ \mc{U}_j \}_{j\in J}$ of families of a set $X$, describe the generalized topology they generate.

\vspace{2mm}
\textbf{Acknowledgements.} This paper was mainly written during my stay at the Fields Institute during the Thematic Program on o-minimal Structures and Real Analytic Geometry in 2009.
I thank the Fields Institute for their warm hospitality.

{\sc Politechnika Krakowska

Instytut Matematyki

Warszawska 24

PL-31-155 Krak\'ow

Poland}

E-mail: \textit{pupiekos@cyf-kr.edu.pl}
\end{document}